\documentclass[a4paper,12pt,legno]{article}

\usepackage{amsmath}
\usepackage{amsfonts}
\usepackage{amsthm}
\usepackage{mathrsfs}
\usepackage[T1]{fontenc}
\usepackage{amsmath}
\usepackage{amsfonts}
\usepackage{amsthm}
\usepackage{mathrsfs}
\usepackage{amssymb}

\newtheorem{thm}{Theorem}
\newtheorem{prop}[]{Proposition}

\newtheorem{cor}{Corollary}

\newtheorem{rem}[]{Remark}

\newcommand{\Rset}{\mathbb{R}}
\newcommand{\Cset}{\mathbb{C}}

\usepackage[T1]{fontenc}
\usepackage{amsmath}
\usepackage{amsfonts}
\usepackage{amsthm}
\usepackage{mathrsfs}
\usepackage{amssymb}

\begin{document}

\title{On a method of introducing  free-infinitely divisible probability measures\footnote{Research funded by Narodowe Centrum Nauki (NCN), grant no Dec2011/01/B/ST1/01257.}}

\author{Zbigniew J. Jurek  (University of Wroclaw)}

\date{December 14, 2014}

\maketitle

\begin{quote}
\textbf{Abstract.} Random integral mappings $I^{h,r}_{(a,b]}$ give isomorphisms between the sub-semigroups of the classical  $(ID, \ast)$ and the free-infinite divisible $(ID,\boxplus)$  probability measures.  This  allows us to introduce new  examples of  such  measures  and their corresponding characteristic functionals.



\emph{Abbreviated title: Free-infinitely divisible measures }

\end{quote}

In this paper we prove Fourier type characterizations for  new classes of
limiting distributions (subsets of infinitely divisible  laws) in the free probability theory.
One of them (Proposition 1) is the free  analog of the class $\mathcal{U}$. The class $\mathcal{U}$ was
defined as the class of weak  limits for sequences of the form:
\[
U_{r_n}(X_1)+U_{r_n}(X_2)+...+U_{r_n}(X_n)+x_n \Rightarrow \mu ,  \ \ \qquad \qquad \ (\star)
\]
where random variables $(X_n)$ are stochastically independent, the above summands are infinitesimal,
$(x_n)$ are real numbers and \emph{the non-linear shrinking deformations} $U_r, r>0$, are defined as follows;
\begin{equation*}
U_r(0):=0 ,  \ \ \ U_r(x):=\max\{|x|-r,0\}\,\frac{x}{|x|} , \ \mbox{for} \ \ x\neq 0.
\end{equation*}
Probability measures $\mu$ in $(\star)$ are called \emph{s-selfdecomosable} ( $"s"$, because of the shrinking operators
$U_r$.)
They were introduced for measures on Hilbert spaces in Jurek (1977) and in Jurek (1981). Later on studied in Jurek (1984) and (1985). More recently, in Bradley and Jurek (2014), the Gaussian limit
in $(\star)$ was proved in a case when the independence  was replaced by some strong mixing conditions. On other hand, in Arzimendi and Hasebe (2014), Section 5,  measures in class $\mathcal{U}$ have studied via the unimodality property of their L\'evy spectral measures.

Replacing in $(\star)$ the $U_r's$ by the linear dilations $T_r(x);=r\,x$ we get the L\'evy class $L$ of so called \emph{selfdecomposable} distributions.
In particular,  we obtain \emph{stable  distributions}, when $X_i's$ are also identically  distributed.

For purposes in this paper we need the descriptions of the classes $L$ and $\mathcal{U}$ in terms
the random integral representations  like the ones in (4);
cf.  Jurek and Vervaat (1983) and Theorem 2.1 in Jurek (1984), respectively.

\medskip
\textbf{0. The isomorphism.} \ \ Traditionally, let $\phi_{\mu}$ denotes the Fourier transform (the characteristic function)
of a probability measure $\mu$   and let   $V_{\nu}$ denotes the Voiculescu transform
of a probability measure $\nu$ (the definition is given in the subsection  \textbf{1.2.} below). Then
for a class $\mathcal{C}$ of  classical $\ast$ - infinitely divisible probability measures we define
its free $\boxplus$ -infinitely divisible counter part $\mathfrak{C}$ as follows:
\begin{equation}
\mathfrak{C}=\{\nu: V_{\nu}(it)=it^2\int_0^{\infty}\log\phi_{\mu}(-v)e^{-tv}dv, \ t>0;
\ \ \mbox{for some}\ \mu\in\mathcal{C}\}
\end{equation}
Conversely, for a class $\mathfrak{C}$ of  $\boxplus$ -infinitely divisible  measures we
 define its classical $\ast$ - infinitely divisible
counterpart $\mathcal{C}$ as follows:
\begin{equation}
\mathcal{C}=\{\mu : it^2\int_0^{\infty}\log\phi_{\mu}(-v)e^{-tv}dv= V_{\nu}(t), \ t>0;
\ \ \mbox{for some}\ \nu\in\mathfrak{C}\}
\end{equation}
It is notably that above, and later on, we consider $V_{\nu}$  (and  Cauchy transforms)
only on the  imaginary axis. Still, it is sufficient to perform the explicit  inverse procedures; cf. Section \textbf{1.3.} below.

We illustrate the relation between  classes $\mathcal{C}$ and $\mathfrak{C}$ (in fact, an isomorphism)
via examples and will prove among others that :
\begin{multline*}
\nu \ \ \mbox{\emph{is $\boxplus$ s-selfdecomposable  if and only if for t>0}}\\
 V_{\nu}(it)=\frac{a}{2}+\frac{\sigma^2}{3}\frac{1}{it}+
\int_{\Rset\setminus{\{0\}}}\Big(\frac{(it)^2\,[\log(it-x)-\log(it)]}{x}-it-\frac{1}{2}x1_B(x)\Big)M(dx)
\end{multline*}
for some constants $a\in\Rset, \sigma^2\ge0$ (variance) and  a L\'evy spectral measure $M$. The parameters $a, \sigma^2$ and the measure $M$  correspond to $\mu=[a,\sigma^2, M]$ in (1) and (2);
for other details cf. Proposition 1.

\medskip
The method (idea) of
inserting the same characteristics ( parameters) into different integral  kernels
can be traced to Jurek and Vervaat (1983), p. 254, where it was used to describe the class L of selfdecomposable distributions.
Namely, $[a, R, M]_L$ (the triple in the L\'evy-Khintchine formula)  was identified with $[a,R,M]$. For other
such examples  cf. Jurek (2011),
an invited talk at 10th Vilnius Conference. Similarly, Bercovici  and Pata  (1999)
introduced a bijection between the semigropus of classical and free
infinite divisible probability measures. In this paper the identification
is done on the level of L\'evy exponents $\log \phi_{\mu}$ (cummulants) and
Voiculescu transforms $V_{\nu}$.

In this paper we show the explicit relation between $V_{\textbf{w}}(z)=1/z$ and $\phi_{N(0,1)}(t)=\exp(-t^2/2)$, where $\textbf{w}$ is the Wigner semicircle distribution and $N(0,1)$ is the standard Gaussian measure.

\medskip
\medskip
\textbf{1. The classical $\ast$- and the $\boxplus$- free infinite divisibility.}

\medskip
\textbf{1.1.}\  A probability measure $\mu$ is \emph{ $\ast$-infinitely divisible} (ID, $\ast$)  if for each natural $n\ge2$ there  exists  probability measure $\mu_n$ such that $\mu_n^{\ast n}=\mu$. Equivalently, its characteristic function $\phi_{\mu}$ (Fourier transform) admits the following form (L\'evy-Khinchine formula)
\begin{equation}
\log\phi_{\mu}(t)= ita-\frac{1}{2}\sigma^2t^2+\int_{\Rset\setminus\{0\}}\big(e^{itx}-1-itx 1_{|x| \leq 1}\big)M(dx),  \ \ t \in \Rset,
\end{equation}
and the triplet $a\in\Rset,\sigma^2\ge0$ (covariance)  and a positive Borel measure $M$ are uniquely determined by $\mu$; in short we write $\mu=[a,\sigma^2,M]$. A sigma-finite measure $M$ in (3)  is finite on all open complements of zero and integrates
$|x|^2$ in every finite neighborhood of zero. It is called  \emph{ the  L\'evy spectral measure} of $\mu$.

In recent year has been considerable interest in studying \emph{random integral representations} of infinitely divisible probability measures or their L\'evy measures $M$; cf. Jurek (2012) and  references therein. Namely,
for a continuous $h$ and a monotone  right continuous $r$,on an interval $(a,b]$,  one defines
\begin{equation}
I^{h,r}_{(a,b]}(\nu):= \mathcal{L}\big(\int_{(a,b]}h(s)dY_{\nu}(r(s)\big),
\end{equation}
where $\mathcal{L}(X)$ denotes the probability distribution of $X$  and $Y_{\nu}$ is a cadlag L\'evy process such that $\mathcal{L}(Y_{\nu}(1))=\nu$.

In terms of characteristic  functions (4) means that\begin{equation}
\widehat{(I^{h,r}_{(a,b]}(\nu))}(t)=\exp\int_{(a,b]}\log\phi_{\nu}(\pm\,h(s)\,t)(\pm)dr(t), \ \ t \in\Rset,
\end{equation}
where the minus sign is for decreasing $r$  and plus for increasing $r$; cf. Jurek and Vervaat (1983), Lemma 1.1 or  Jurek (2007) (in the proof of Theorem 1) or Jurek (2012). Moreover,  $\widehat{(I^{h,r}_{(a,b]}(\nu))}$ denotes here the characteristic function  of the probability measure $I^{h,r}_{(a,b]}(\nu)$. In (5), we may write $\phi_{\nu}(-w)=\phi_{\nu^{-}}(w), w\in\Rset$, where $\nu^{-}$ is the reflected measure, that is, $\nu^{-}(B):=\nu(-B)$ for all Borel sets $B$.

\medskip
For the purposes below we consider the following specific random integral mapping:
\begin{equation}
(ID, \ast)\ni \mu\to \mathcal{K}(\mu)\equiv I^{s,1-e^{-s}}_{(0,\infty)}(\mu)=
 \mathcal{L}\big(\int_0^{\infty}\,s\,dY_{\mu}(1-e^{-s})\big)\in\mathcal{E},
\end{equation}
from Jurek (2007), formula (17). There, it was done for any real separable Hilbert space.
(In Barndorff-Nielsen and Thorbjornsen (2006), and in other works,  the mapping (6) was denoted by the letter $\Upsilon$
and (originaly)  was defined on the family of L\'evy measures on a real line).

The mapping $\mathcal{K}$ is an isomorphism between convolution semigroup $ID$ and $\mathcal{E}$ (range of the mapping $\mathcal{K}$). Moreover, if
$\mu=[a, \sigma^2, M]$ then from (3) we get
\begin{equation}
\phi_{\mathcal{K}(\mu)}(t)=\exp\{ita-\sigma^2 t^2 +
\int_{\Rset\setminus\{0\}}\big( \frac{1}{1-itx}-1- itx1_{\{|x|\le 1\}}(x \big)M(dx)\};
\end{equation}
cf. Jurek (2006), Corollary 5. For a general theory of the calculus on random integral mappings of the form (2)  cf. Jurek (2012).

\medskip
\textbf{1.2.} \ D. Voiculescu and others studying so called
\emph{free-probability} introduced new binary operations on
probability measures and termed them accordingly \emph{free-convolutions}; cf.
Bercovici-Voiculescu (1993) and references therein. To recall the definition $\boxplus$ convolution we need some auxiliary notions.

For a measure $\nu$, its \emph{Cauchy transform} is given as follows
\begin{equation}
G_{\nu}(z):=\int_{\Rset}\frac{1}{z-x}\nu(dx), \ \   \  z \in \Cset\setminus \Rset.
\end{equation}
Furthermore, having $G_{\nu}$ we define $F_{\nu}(z):=1/G_{\nu}(z)$ and  then \emph{the Voiculescu transform} as
\begin{equation}
V_{\nu}(z):= F^{-1}_{\nu}(z)-z, \ z\in \Gamma_{\eta,M}:=\{x+iy\in \Cset^+: |x|<\eta y, y>M\},
\end{equation}
where a such region (called Stolz angle) exists and the inverse function is well defined on it; cf. Bercovici and Voiculescu (1993), Proposition 5.4 and Corollary 5.5.

The functional $\nu\to V_{\nu}(z)$ is an analogue  of the classical Fourier transform $\nu\to \phi_{\nu}(t), t\in\Rset$.
The fundamental  fact is that, for two measures $\nu_1$ and $\nu_2$ one has
$$V_{\nu_1}(z)+V_{\nu_2}(z)=V_{\nu_1
\boxplus \nu_2}(z),$$
for a uniquely determined probability measure, denoted as
$\nu_1\boxplus \nu_2$.
This property allowed to
introduce the notion of  $\boxplus$ free-infinite divisibility.
For this new $\boxplus$ infinite divisibility we have  the following analog of the  L\'evy-Khintchine formula (3):
\begin{thm}
A measure $\nu$ is $\boxplus$ infinitely divisible if and only if
\begin{equation}
V_{\nu}(z)= b + \int_{\Rset}\,\frac{1+sz}{z-s}\,\rho(ds), \qquad z\in\Cset\setminus \Rset,
\end{equation}
for some uniquely determined real constant $b$ and a finite Borel measure $\rho$.
\end{thm}
Cf. Bercovici and Voiculescu (1993). The integral formula (10), in complex analysis, is
called \emph{the Nevanlinna-Pick} formula.
\begin{rem}
\emph{Note that $V_{\nu} : \Cset^+\to \Cset^-$ is an analytic function and
for the mappings   $T_c\,x:=c\,x, x\in \Rset$,  $(c>0)$ and the image measure $T_c\nu$  we have
$V_{T_c\nu} (z)=c\,V_{\nu}(c^{-1}z)$, for z in appropriate  Stolz angle;
cf. Bercovici-Voiculescu (1993) or Barndorff-Nielsen (2006), Lemma 4.20. This is in contrast to characteristic functions of measures where we have
$\phi_{T_c\mu}(t)=\phi_{\mu}(ct)$, for all $t\in\Rset$.}
\end{rem}

\medskip
\textbf{1.3.} \ \
For some analogies and comparison below, let us recall from Jurek (2006) that the restricted versions 
of  $G_{\nu}$ and $V_{\nu}$ are just those functions considered only on the imaginary axis. Then we have  that
\begin{equation}
\frac{1}{it}G_{\nu}(\frac{1}{it})=\int_{\Rset}\frac{1}{1-itx}\,\nu(dx)
=\phi_{e\cdot  \eta}(t),  \  \ t \neq 0, \ \mbox{and} \  \lim_{t\to 0}\frac{1}{it}G_{\nu}(\frac{1}{it})=1
\end{equation}
where $ e\cdot \eta$ means  the probability distribution of a product of stochastically independent
 rv's: the standard exponential $e$ and the variable $\eta$ with probability distribution $\nu$.

The identity (11) means that we can retrieve a measure $\nu$ from the characteristic function $\phi_{e\cdot \eta}$;
cf. Jurek (2006), the proof of Theorem 1, on p.189 and Examples on pp. 195-198.
This is in sharp contrast with the classical Stieltjes inversion formula where one needs
to know $G_{\nu}$ in strips  of complex plane; see it, for instance, in Nica and Speicher (2006), p. 31.

In fact, we have even more straightforward relation. Namely,
\begin{equation}
\int_0^{\infty}\phi_{\nu}(s)e^{-ts}ds=\overline{iG_{\nu}(it)}, \ \ \mbox{for} \ \ t>0,
\end{equation}
cf. Jankowski and Jurek (2012), Proposition 1.
Thus, restricted Cauchy transforms
are just Laplace transforms of characteristic functions.

In the spirit of (11) and (12),  instead of (10), let us introduce  \emph{the restricted Voiculescu transform} as
\begin{equation}
V_{\nu}(it)\equiv k_{b,\rho}(it):= b +\int_{\Rset}\frac{1+it\,s} {it-s}\,\rho(ds), \ \  t\neq 0.
\end{equation}
From  the inversion formula in Theorem 1 in Jankowski and Jurek (2012) and from (13),  we have that
\begin{multline}
b= \Re k_{b,\rho}(i);  \qquad \rho(\Rset):= - \Im k_{b,\rho}(i);   \ \ \ \mbox{and for the measure $\rho$ we have}  \\
\int_0^{\infty}\phi_{\rho}(r)e^{-w\,r}dr=
\frac{ i\, k_{b,\rho}(-iw) - i\,\Re k_{b,\rho}(i) - w\Im k_{b,\rho}(i)}{w^2-1},   \ \ \  w>0.
\end{multline}
(See there also the comment (a paradigm) in the first paragraph in the introduction on p. 298.)

In order to have (13) in a form more explicitly related to (7),
let us define the new triple: shift a, the variance $\sigma^2$ and the L\'evy spectral measure $M$, as follows:
\begin{multline}
 \sigma^2:=\rho(\{0\}) ; \ \  \qquad \  \  \rho(dx):=\frac{x^2}{1+x^2} M(dx), \  \
\mbox {on} \   \Rset\setminus{\{0\}};  \\
a:= b+\int_{\Rset} s\ \big(1_{|x|\le 1}(s)-\frac{1}{1+s^2}\big)\ \rho(ds); \qquad \qquad \qquad  \qquad \qquad \qquad
\end{multline}
Then from (13), with some calculations, we get
\begin{equation}
it\,k_{b,\nu}(\frac{1}{it})=
ia\,t -\sigma^2 t^2+ \int_{\Rset}\big(\frac{1}{1-sit}-1-its\,1_{|s|\le 1}\big)\,M(ds), \ \mbox{for} \ t <0;
\end{equation}
For computational details cf. Barndorff-Nielsen and Thorbjornsen (2006), Proposition 4.16, p. 105

\medskip
\textbf{1.4.} \ The functions $G_{\nu}(z)$ and $V_{\nu}(z)$ are analytic in some complex domains and thus
are uniquely determined by their values on imaginary axis  (more generally, on subsets with
limiting points in their domains).

\noindent [Note that for $p(z):= i \Im z$ and $q(z):=z$ we have  that $p(it)=q(it)$ although they are different.
Of course, $p$ is not an analytic function! ]

\medskip
\textbf{1.5.} \ \ Because of (6), (7) and (16),  here is the explicit  relation (an isomorphism) between
the free- $\boxplus$  and the classical - $\ast$ infinite divisibility:
\begin{thm}
A probability measure $\nu$ is $\boxplus$-infinitely divisible if and
only if there exist a unique $\ast$-infinitely divisible probability
measure $\mu$ such that
\begin{multline}
(it)\,V_{\nu}((it)^{-1})=\log\big( I^{s,\,
1-e^{-s}}_{(0,\infty)}(\mu)\big)^{\widehat{}}(t)=\\ \log
\Big(\mathcal{L}
(\int_0^{\infty}s\,dY_{\mu}(1-e^{-s}))\Big)^{\widehat{}}(t)
=\int_0^{\infty}\log\phi_{\mu}(ts)e^{-s}ds ,
\  \  \mbox{for} \ \  t<0,
\end{multline}
where $(Y_{\mu}(u), u\ge 0)$ is a L\'evy process such that
$\mathcal{L}(Y_{\mu}(1))=\mu$.

Equivalently, we have that for $\boxplus$-infinitely divisible $\nu$ its Voiculescu transform $V_{\nu}$ is of the form
\begin{equation}
V_{\nu}(it)=it\int_0^{\infty}\log \phi_{\mu}(-t^{-1}s)e^{-s}ds=it^2\int_0^{\infty}\log\phi_{\mu}(-v)e^{-tv}dv, \ t>0;
\end{equation}
for an uniquely determined $\ast$- infinitely divisible measure $\mu$.
\end{thm}
This is a rephrased  version of   Corollary 6 in Jurek (2007).

\medskip
Statements  (17) and (18) are equivalent as one can be deducted from the other.
Moreover, they  provide easy way of computing examples (classes) of free
$\boxplus$-infinitely divisible measures  from their counterparts in $(ID,\ast)$;
cf. Propositions  1- 4 below.

Also note that in both cases (17) and (18) we have  Laplace transform of functions
$\log\phi_{\mu}(- t)$ of $\ast$-infinitely divisible measures $\mu$.
\begin{rem}
\emph{ From the  first line in (17) we see that measures from $(ID,\boxplus)$
can be identified  with measures from  the  semigroup $\mathcal{E}=I^{s,1-e^{-s}}_{(0,\infty)}(ID)$.}
\end{rem}

\medskip
\textbf{2. Examples of explicit  relations between free $\boxplus$ - and classical $\ast$ - infinite divisible probability measures.}

\medskip
In the first three subsections, for a given class $\mathcal{C}$ classical $\ast$-infinitely
divisible measures we identify its counterpart  $\mathfrak{C}$ of free
$\boxplus$ - infinitely divisible companions.

\medskip
\textbf{2.1.} For the free $\boxplus$ analog of s-selfdecomposable distributions we have
\begin{prop}
A probability distribution $\mathfrak{u}$ is free $\boxplus$ \emph{s-selfdecomposable},
in symbols, $\mathfrak{u} \in (\mathcal{U}, \boxplus)$,  if and only it
there exist a unique $\mu=[a,\sigma^2,M] \in (ID, \ast)$ such that its Voilculescu transforms
have representations
\begin{equation}
(a) \ \ it\,V_{\mathfrak{u}}( \frac{1}{it})=
\log \widehat{(I^{v,r_u(v)}_{(0,\infty)}(\mu^{-}))}(t), \  \ \mbox{with} \ r_u(v):= e^{-v}-v\Gamma(0,v);
\end{equation}
\begin{equation*}
(b)  V_{\mathfrak{u}}(it)=\frac{a}{2}+\frac{\sigma^2}{3}\frac{1}{it}+
\int_{\Rset\setminus{\{0\}}}\Big(\frac{(it)^2\,[\log(it-x)-\log(it)]}{x}-it-\frac{1}{2}x1_B(x)\Big)M(dx)
\end{equation*}
for $t>0$.

(c)\ For $z\in\Cset^+$ we have
\begin{equation*}
V_{\mathfrak{u}}(z)=\frac{a}{2}+\frac{\sigma^2}{3}\frac{1}{z}+
\int_{\Rset\setminus{\{0\}}}\Big(z^2\frac{[\log(z-x)-\log(z)]}{x}-z -\frac{1}{2}x1_B(x)\Big)M(dx).
\end{equation*}
\mbox{Above} $\Gamma(0;x):=\int_x^{\infty}\frac{e^{-s}}{s}ds, x>0,$ \mbox{is the incomplete Euler gamma function.}
\end{prop}

\medskip
\emph{Proof.}
Recall that $\lambda$ is  classical $\ast$-s-selfdecomposable, i.e., $\lambda\in (\mathcal{U}, \ast)$ if and only if
$\lambda=I^{s,\,s}_{(0,1]}(\mu)$ for some $\mu\in ID$; cf. Jurek (1984), Theorem 2.1
for different characterizations of this class or Jurek (1985).  Thus for $w\neq 0$, using (3), we get
\begin{multline}
\log\phi_{\lambda}(w)=\log\big( I^{s,\,s}_{(0, 1]}(\mu)\big)^{\widehat{}}(t)
=\int_0^1\log\phi_{\mu}(w\,u)du\\=i\frac{1}{2}aw-\frac{1}{6}\,\sigma^2 w^2+
\int_{\Rset\setminus\{0\}}[\frac{e^{iwu}-1}{iwu}-1 -i \frac{1}{2}wu\,1_{|u|\le1}(u)]\,M(du).
\end{multline}
[The formula (20) is as in Corollary 7.1 in Jurek (1984). However, it was obtained there by using Choquet's Theorem on extreme points in a subset of all L\'evy spectral measures.]

From  (17) and the first line in (20)  we get
\begin{multline}
it\,V_{\mathfrak{m}}(\frac{1}{it})=\int_0^{\infty}\log\phi_{\lambda}(t\,s)e^{-s}ds=
 \int_0^{\infty}\int_0^{1}\log\phi_{\mu}(tus)e^{-s}duds,  \ \ (v:=su)\\
=\int_0^{\infty}\int_0^s\log \phi_{\mu} (tv)\frac{e^{-s}}{s}dv\,ds=\int_0^{\infty}\log\phi_{\mu}(tv)\Gamma (0,v)dv.
\end{multline}
Let us define the (decreasing) time change $r_u(v)$ for $v>0$ as follows
\begin{multline*}
r_u(v):=\int_v^{\infty}\Gamma(0,w)dw =\int_v^{\infty}\int_w^{\infty}\frac{e^{-s}}{s}ds\,dw \\ = \int_v^{\infty}\int_v^{s}\frac{e^{-s}}{s}dw\,ds =\int_v^{\infty} \frac{e^{-s}}{s}(s-v)ds=e^{-v}-v\Gamma(0,v).
\end{multline*}
Then
taking into account (5) and putting $r_u$  into (21) to get
\begin{equation*}
it\,V_{\mathfrak{m}}(\frac{1}{it})=\int_0^{\infty}\log\phi_{\mu^{-}}(-tv)(1)dr_u(v) =\log(\widehat{I^{v,\,r_u(v)}_{(0,\infty)}(\mu^{-})})(t)
\end{equation*}
and this completes the proof of the part (a).

For part (b),  using  (17), (3)  and the first line in (20), after interchanging the order of integration,  we get for $t<0$,
\begin{multline}
it\,V_{\mathfrak{m}}(\frac{1}{it})=\int_0^{\infty}\log\phi_{\lambda}(t\,s)e^{-s}ds= \int_0^1\int_0^{\infty}\log\phi_{\mu}(tus)e^{-s}dsdu \\  =
\frac{1}{2}
iat-\frac{1}{3}\sigma^2t^2+ \int_{\Rset\setminus\{0\}}\int_0^1[\int_0^{\infty}\big(e^{itusx}-1-itusx1_B(x)\big)e^{-s}ds]du M(dx)\\=
\frac{1}{2}
iat-\frac{1}{3}\sigma^2t^2+ \int_{\Rset\setminus\{0\}}[\int_0^1\frac{itux}{1-itux}du -\frac{1}{2}itux1_B(x)]M(dx)
\\
=\frac{1}{2}
iat-\frac{1}{3}\sigma^2t^2+ \int_{\Rset\setminus\{0\}}[ - \frac{\log(1-itx)}{itx}-1-\frac{1}{2}itx1_B(x)]M(dx).
\end{multline}
Substituting $-1/t$ for $t$ in (22) we arrive at
\[
V(it)=\frac{a}{2}+\frac{\sigma^2}{3}\frac{1}{it}+
\int_{\Rset\setminus{\{0\}}}\Big((it)^2\frac{[\log(it-x)-\log(it)]}{x}-it-\frac{1}{2}x1_B(x)\Big)M(dx),
\]
which gives (b).
Part (c) is an analytic extension of (b) and this completes a proof of Proposition 1.

\medskip
\noindent[Equality (b) can be also obtained by putting second line in (20) into (18).]

\begin{rem}
\emph{ (i) There is  the connection between the class $(\mathcal{U},\boxplus)$ and the class
$I^{v,r_u(v)}_{(0,\infty)}(ID)$ via (a) in Proposition 1.}

\emph{(ii) $I^{t, 1-e^{-t}}_{(0,\infty)}\circ I^{s,s}_{(0,1]} =I^{v,r_u(v)}_{(0,\infty)} \ \mbox{with} \ r_u(v):= e^{-v}-v\Gamma(0,v);$}
\end{rem}

\medskip
Let $h_1\otimes h_2$ denotes the tensor product of functions and let $\rho_1\times \rho_2$
denotes the product of measures. Then we have
\begin{cor}
For $h_1(t):= t\,1_{(0,1)}(t), \rho_1(dx):= 1_{(0,1)}(x)dx$ and $h_2(s):=s\,1_{(0,\infty)}(s)$, and  $\rho_2(dy):=e^{-y}dy$ we have that
\[
(h_1\otimes h_2)(\rho_1\times \rho_2)(dw)= 1_{(0,\infty)}(w)\Gamma (0,w)dw
\]
\end{cor}
This is in fact the calculation performed in (21) starting with the first double integral.

\medskip
\textbf{2.2.} For free $\boxplus$ - selfdecomposability we have the following:
\begin{prop}
A probability distribution $\mathfrak{s}$ is $\boxplus$-\emph{selfdecomposable}, in symbols, $\mathfrak{s}\in (L,\boxplus)$,  if and only it there exist a unique
$\mu=[a,\sigma^2,M] \in (ID_{\log}, \ast)$ such that
\begin{equation}
(a) \ itV_{\mathfrak{s}}( \frac{1}{it})= \log \widehat{\big(I^{w, \Gamma(0,w)}_{(0,\infty)}(\mu^{-})\big)}(t), \ \mbox{for}
 \ t<0; \ \ \Gamma(0, w):=\int_w^{\infty}\frac{e^{-s}}{s}\,ds
\end{equation}
Equivalently,
\begin{equation}
(b)  \ \ V_{\mathfrak{s}}(it)= a +\frac{\sigma^2}{2}\frac{1}{it} +\int_{\Rset\setminus\{0\}}[it \ln\frac{it}{it-x}-x1_{\{|x|\le 1\}}]M(dx),  \ \  \mbox{for}  \  \ t>0.
\end{equation}
That is, for $z\in\Cset^+$,
\[
(c)  \ \ \ V_{\mathfrak{s}}(z)= a +\frac{\sigma^2}{2}\frac{1}{z} +\int_{\Rset\setminus\{0\}}[z \ln\frac{z}{z-x}-
x1_{\{|x|\le 1\}}]M(dx).
\]
\end{prop}

\medskip
\emph{Proof.} Recall that $\rho\in (L, \ast)$, in other words, $\rho$ is $\ast$- selfdecomposable if and only if
$\rho=I^{e^{-s}\,,s}_{(0,\infty)}(\mu)$ for some $\mu\in ID_{\log}$; cf. Jurek-Vervaat (1983) or Jurek-Mason (1993),
Chapter 3 (with operator $Q=I$). Hence
\begin{multline}
\log\phi_{\rho}(w)=
\log
\big( I^{e^{-s},\,s}_{(0,\infty)}(\mu)\big)^{\widehat{}}(w)=
 \int_0^{\infty}\log\phi_{\mu}(we^{-u})du
\\=iaw-\frac{1}{4}\,\sigma^2 w^2+
\int_{\Rset\setminus\{0\}}\Big(\int_{(0,1)}\frac{e^{iwrx}-1}{r}\,dr-iwx\,1_{|x|\le1}(x)\Big)\,M(dx).
\end{multline}
[The characterization in (25) of selfdecomposable distributions was first obtained by K. Urbanik
(by the method of extreme points) and then by Jurek and Vervaat (1983), formula (4.5) on p. 255.]

From (17) and the first line in (25), for $t<0$,  we have
\begin{multline}
itV_{\mathfrak{s}}( \frac{1}{it})=
\int_0^{\infty}\log\phi_{\rho}(ts)e^{-s}ds=\int_0^{\infty}\int_0^{\infty}\log\phi_{\mu}(tse^{-u})e^{-s}duds  \ (w:=se^{-u} ) \\
\int_0^{\infty}[\int_0^s\log\phi (tw)\frac{1}{w}dw]e^{-s}ds=\int_0^{\infty}\log \phi(tw)\frac{1}{w}
[\int_w^{\infty}e^{-s}ds]dw\\
=\int_0^{\infty}\log\phi_{\mu}(tw)\frac{e^{-w}}{w}dw=
\int_0^{\infty}\log\phi_{\mu^{-}}(-tw)(-1)(\Gamma(0,w))^{\prime}dw\\
= \log {\widehat{ \big(I^{w, \ \Gamma(0,w)}_{(0,\infty)}(\mu^{-})\big)}}(t), \ \ \ \ \ \ \ \ \
\end{multline}
(see the equality (5), which gives (a).

Or equivalently, using the second line in (25) and (3) we get
\begin{multline*}
itV_{\mathfrak{s}}( \frac{1}{it})  =\int_0^{\infty}\Big[iats-\frac{1}{4}\,\sigma^2 (ts)^2 \\ +
\int_{\Rset\setminus\{0\}}\Big(\int_{(0,1)}\frac{e^{itsrx}-1}{r}\,dr-itsx\,1_{|x|\le1}(x)\Big)\,M(dx)\Big]e^{-s}ds \\
=iat-\frac{1}{2}\sigma^2t^2 +
\int_{\Rset\setminus\{0\}}\Big(\int_{(0,1)}(\frac{1}{1-itxr}-1)\,\frac{1}{r}\,dr-itx\,1_{|x|\le1}(x)\Big)\,M(dx)\\
=iat-\frac{1}{2}\sigma^2t^2 +
\int_{\Rset\setminus\{0\}}\Big(\int_{(0,1)}\frac{itx}{1-itxr}dr-
 \,itx\,1_{|x|\le1}(x)\Big)M(dx)  \\
= iat-\frac{1}{2}\sigma^2t^2 +
\int_{\Rset\setminus\{0\}}\Big(-\ln(1-itx)-
 \,itx\,1_{|x|\le1}(x)\Big)M(dx).
\end{multline*}
Now substituting for $s:= -1/t>0$  w get
\begin{equation*}
\frac{1}{is}V_{\mathfrak{s}}(is)=\frac{-ia}{s}-\frac{\sigma^2}{2s^2}+\int_{\Rset\setminus\{0\}}\Big(-\ln(1+ix/s)+
 \,ix/s\,1_{|x|\le1}(x)\Big)M(dx),
\end{equation*}
and hence the equality (b).

Part (c) is an analytic continuation of (b) and this concludes a proof of Proposition 2.

In terms of tensor product we have
\begin{cor}
For $h_1(t):=t,  \ \rho_1(dt):=e^{-t}, \ 0<t<\infty$ and $h_2(s):=e^{-s}; \ \rho_2(ds):=ds , 0<s<\infty$ then
\[
(h_1\otimes h_2)(\rho_1 \times \rho_2)(du)= \frac{e^{-u}}{u}du
\]
\end{cor}
This follows from (26) starting from the double integral and the calculations that follows. See also Jurek (2012).

\begin{rem}
\emph{(a) From (26) we get  $I^{t,1-e^{-t}}_{(0,\infty)} \circ I^{e^{-t},t}_{(0\infty)}(\mu)=
I^{w,\ \Gamma(w)}_{(0,\infty)}(\mu)$, for $\mu \in ID_{\log}$. This is Thorin class $\mathcal{T}$ and  thus $\mathcal{T}$ can be identified with the class of free-selfdecomposable measures. (Compare Remarks 2 and 3(i).) }

\emph{(b) For  \emph{the standard exponential measure} $\mathbf{e}$ and \emph{ standard Poisson measure $e(\delta_1)$}  we have
\begin{equation}
I^{w,\ \Gamma(w)}_{(0,\infty)}e(\delta_1))= \mathbf{e}.
\end{equation}}
\end{rem}
To see (b) note that
\[
\log \phi_{\mathbf{e}}(v)=\int_0^{\infty}\,(e^{ivx}-1)\,\frac{e^{-x}}{x}dx, \ \  \  \mbox{for} \ v\in\Rset,
\]
therefore for $\emph{e}(\delta_1)$ we conclude
\begin{multline*}
\log {\widehat{\big(I^{w, \ \Gamma(w)}_{(0,\infty)}(e(\delta_1)\big)}}(t) \\
=\int_0^{\infty}\log\phi_{e(\delta_1)}(t\,w)\frac{e^{-w}}{w}dw
=\int_0^{\infty}\,(e^{itw}-1)\,\frac{e^{-w}}{w}dw = \log \phi_{\mathbf{e}}(t),
\end{multline*}
which proves that (unexpected ?) relation (27) between Poisson (discrete) and exponential (continuous) distributions.

\medskip
\medskip
\textbf{2.3.} For free $\boxplus$-stable distributions we have:
\begin{prop}
A measure $\nu$ is non-Gaussian free-$\boxplus$ stable  if and only if for $t>0$ its
Voiculescu transform $V_{\nu}$ should be such that it is EITHER
\begin{equation}
V_{\nu}(it)=
a- C\frac{\Gamma{(2-p)}}{1-p}\Gamma{(p+1)}\cos \frac{\pi p}{2}\ \big[ i^{p-1}(it)^{1-p}
(i-\beta \tan(\frac{\pi p}{2})\big],
\end{equation}
where $ a\in \Rset, C>0$, $0<p<1$ or $1<p<2$ and $|\beta|\le1$,

OR $p=1$ and
\begin{equation}
V_{\nu}(it)=a-C\beta (1-\gamma) +\frac{C}{2}\big[2 \beta\log (it)-i\pi(1+\beta)  \big]
\end{equation}
where $1-\gamma= \int_0^{\infty}w \log w \, e^{-w}dw$ (Euler constant $\gamma \sim 0.577$).
\end{prop}

\emph{Proof.}
For classical $\ast$-stable measures from Meerscheart and Scheffler (2001), Theorem 7.3.5,  p. 265  we have that
\begin{multline}
\mu \ \ \mbox{\emph{is non-Gaussian  $\ast- $stable iff and only if there exist}} \\  C>0, \
a\in\Rset, \ 0<p<1,\, 1<p<2, \ -1 \le \beta\le 1 \  \mbox{\emph{such that for each}} \  t\in\Rset \qquad \qquad \qquad  \\
\log\phi_{\mu}(t)= ita-C\,\frac{\Gamma(2-p)}{1-p}\, \cos(\frac{\pi p}{2})\,|t|^p
\Big(1 -i \beta \,sign (t) \tan(\frac{\pi p}{2})\Big);
\end{multline}
and for $p=1$ we have
\begin{equation}
\log\phi_{\mu}(t)= ita-C\,\frac{\pi }{2}\,|t|\Big( 1+i\,\beta\frac{2}{\pi} sign(t)\log|t|\Big), \  t\in\Rset \,,
\end{equation}
\noindent where $\beta:= 2\theta-1$ is the skewness parameter; $0\le \theta\le 1$ is the probability of the positive tail of L\'evy measure $M$ of $\mu$, that is, for $r>0$, we have $M(x>r)=\theta Cr^{-p}$. And $1-\theta$ is the probability of the negative tail of $M$.

\medskip
In order to get (28) one needs insert (30) into  first equality in (18) and perform some easy calculations. Similarly, putting (31) into (18)
and using the identity $\log i = i \pi/2$ one gets equality (29), which completes a proof of Proposition 3.
\begin{rem}
\emph{(i) \ In some papers and books often there is a small but essential
error. Namely, in (30), there is $\beta$ instead of $(-\beta$); cf. P. Hall (1981).
\newline
(ii) \ Note that the expressions  in square brackets in (28) and (29) are identical, up to the sign, with those in Proposition 5.12 in Bercovici-Pata (1999).
Also compare Biane's formulas for free-stable distributions in the Appendix there.}
\end{rem}

\medskip
\medskip
\textbf{2.4.} \ For a finite Borel measure m, let $e(m):=e^{-m(\Rset)}\sum_{k=0}^{\infty}\frac{m ^{\ast k}}{k!}$ denotes the $\ast$-compound Poison probability measures.
\begin{prop}
A probability measure $\nu$  is $\boxplus$-compound Poison probability measure if and only if
\begin{equation}
V_{\nu}(it)= it \int_{\Rset}\frac{x}{it-x}m(dx) ,  \ \mbox{for} \ t>0,
\end{equation}
for some finite Borel measure $m$ on the real line. Moreover, $\nu$ is free-infinitely divisible if and only if
\begin{equation}
V_{\nu}(z)= b + c^2\,\frac{1}{z}+ \int_{\Rset\setminus{\{0\}}}\,\frac{1}{z-x}m(dx)+z \int_{\Rset} \frac{x}{z-x}m(dx) ,  \ \
z \in \Cset \setminus{\{0\}}
\end{equation}
for some $b, c \in \Rset$ and finite Borel measure $m$ on $\Rset$.
\end{prop}
\emph{Proof.} Since
$\log \phi_{e(m)}(t)= \int_{\Rset} (e^{itx}-1)m(dx), \  t\in\Rset$,
therefore by (17)
\begin{multline*}
V_{\nu}(it)= it \int_0^{\infty}\log\phi_{e(m)}(-t^{-1}s)e^{-s}ds\\=it \int_{\Rset}\int_0^{\infty}(e^{-it^{-1}sx}-1)e^{-s}ds\,m(dx) =it \int_{\Rset}\frac{x}{it-x}m(dx), \ \mbox{for }  \ \ t>0,
\end{multline*}
which completes a proof of (32).
The remaining part is a consequence of above and Theorem 1. [Also see pp. 203-206 in  Nica-Speicher (2006) for the discussion of free compound Poisson distributions]

\medskip
\textbf{2.5.} In this subsection, for a given three examples of $\boxplus$-infinitely divisible measures we identify their classical $\ast$- infinitely divisible companions.

\medskip
\textbf{Example 1.} The probability measure $\textbf{w}$ such that
$V_{\textbf{w}}(z)=\frac{1}{z}, z\neq 0$, is called \emph{free-Gaussian measure}.
Why  such a term?

Note that from Theorem 2, we get
$(it)V_{\textbf{w}}(\frac{1}{it})=-t^2$.
On the other hand, taking standard normal distribution $N(0,1)$ for the measure $\mu$ we get
$$\int_0^{\infty}\log\phi_{\mu}(ts)e^{-s}ds=-t^2\int_0^{\infty}s^2/2e^{-s}ds= - t^2 =
\frac{1}{it}V_{\textbf{w}}(\frac{1}{it}).
$$
So, it is right to call $\textbf{w}$ an analogue of free Gaussian distribution.

More importantly,
$\textbf{w}$ is a weak limit of free-analog of CLT and $\textbf{w}$ is
\emph{ the standard Wigner's semicircle law} with the density $\frac{1}{2\pi}\sqrt{4-x^2}\,1_{[-2,2]}(x)$, mean value zero and variance 1. [Using the inversion formula in (14) for
$V_{\textbf{m}}$ we get $b=0$ and $\rho=\delta_0$ in  (10).]

\medskip
\textbf{Example 2.} The probability measure $\textbf{c}$ with $V_{\textbf{c}}(z)=-i$ is called \emph{free-Cauchy distribution}. Why  such a term?

From Theorem 2, $(it)V_{\textbf{c}}(\frac{1}{it})=t$. On the other hand, taking
the standard Cauchy distribution (with the probability density $\frac{1}{\pi}\frac{1}{1+x^2}$)  for the measure $\mu$ we get
\[
\int_0^{\infty}\log\phi_{\mu}(ts)e^{-s}ds=-|t|\int_0^{\infty}se^{-s}ds=  t=(it)V_{\textbf{c}}(\frac{1}{it}) , \ \mbox{for} \ t<0,
\]
so it justifies the term free - Cauchy measure.
In fact, we have that
\begin{rem}
\emph{ The measure $\textbf{c}$ is the standard Cauchy distribution.
To see that we use the inversion procedure from (14). Thus  $b=0, \textbf{c}(\Rset)=1$ and
\[
\int_0^{\infty}\phi_{\textbf{c}}(r) e^{-wr}dr=\frac{1}{w+1}; \ \ \mbox{i. e.}, \ \phi_{\textbf{c}}(r)=e^{-r}, \ \mbox{for} \ \ r>0.
\]
Consequently, $\phi_{\textbf{c}}(r)=e^{-|r|}$, for  $r\in\Rset$ and hence $\textbf{c}$ is the standard Cauchy probability measure.}
\end{rem}

\medskip
\textbf{Example 3.} The probability measure $\textbf{m}$  such that $V_{\textbf{m}}(z)=\frac{z}{z-1}$ is called \emph{free-Poisson distribution.} Why?

From Theorem 2, $(it)V_{\textbf{m}}(\frac{1}{it})=\frac{it}{1-it}=\frac{1}{1-it}-1$.
On the other hand, if $\mu =e(\delta_1)$ is the standard Poisson distribution  then
$\log\phi_{\mu}(t)=e^{it}-1$,
and by Theorem 2,
$$
\int_0^{\infty}\log\phi_{\mu}(ts)e^{-s}ds=\int_0^{\infty}(e^{its}-1)e^{-s}ds=\frac{1}{1-it}-1= (it)V_{\textbf{m}}(\frac{1}{it}).
$$
In fact, $\textbf{m}$ has the probability density $\frac{1}{2\pi}\sqrt\frac{4-x}{x}
\,1_{(0,4]}(x)$  (so called Marchenko-Pastur law);  cf. Bozejko and Hasebe (2014).

\begin{rem}
\emph{Putting $m:=\delta_1$ in Proposition 4 we retrieve the Example 3.}
\end{rem}

\medskip
\medskip
\textbf{Acknowledgments.} Author would like thank Professor Takahiro Hasebe (Japan) for his help with references to free probability theory.

\bigskip
\bigskip
\noindent {\bf References}

\medskip
O. Arizmendi and T. Hasebe (2014), Classical Scale mixtures of Boolean stable laws, \emph{preprint}.

\medskip
 O.E. Barndorff-Nielsen and S. Thorbjornsen (2006), Classical and Free Infinite divisibility
and L\'evy Processes; \emph{Lect. Notes Math.} \textbf{1866}, pp. 33-159.

\medskip
 H. Bercovici and V. Pata (1999), Stable laws and domains of attraction
in free probability, \emph{Ann. Math.} \textbf{149}, pp. 1023-1060,

\medskip
H. Bercovici and D. Voiculescu (1993), Free convolution of measures with unbounded support,
\emph{Indiana Univ. Math. J.} \textbf{42}, pp. 733-773.

\medskip
M. Bozejko and T. Hasebe (2014),  On free infinite divisibility for classical Meixner
distributions, \emph{preprint}.

 R. C. Bradley and Z. J. Jurek (2014), On central limit theorem for shrunken weakly dependent random variables,  to appear in \emph{Houston J. Math.},  \  arXiv:1410.0214 [math.PR]

\medskip
P. Hall (1981), A comedy of errors: the canonical form for a stable
characteristic function, \emph{Bull. London Math. Soc.} \textbf{13} no 1,
pp. 23-27.

\medskip
L. Jankowski and Z. J. Jurek (2012), Remarks on restricted Nevalinna transforms, \emph{Demostratio Math.} vol. XLV, no2 , pp. 297-307.

\medskip
\noindent Z. J. Jurek (1977),
Limit distributions for sums of
shrunken random variables.
In: \emph{Second Vilnius Conference
on Probability  Theory and  Mathematical Statistics.}
Abstract of Communications 3, pp. 95-96.

\medskip
Z. J. Jurek (1981).
Limit distributions for sums of shrunken
random variables.
{\it Dissertationes Math.\/} vol. \textbf{185},  PWN Warszawa.

\medskip
Z. J.  Jurek (1984),
S-selfdecomposable probability measures as probability distributions of
some random integrals.
{\it Limit Theorems in Probability and
Statistics (Veszpr\'em, 1982)}, Vol.\ I, II,
pp.\ 617-629.
{\it Coll.\ Math.\ Societatis J\'anos Bolyai\/},
Vol. 36.

\medskip
Z. J.  Jurek (1985), Relations between the s- selfdecomposable and selfdecomposable measures, \emph{Ann. Probab.} \textbf{13}, pp. 592-608.

\medskip
Z. J.  Jurek (2006), Cauchy transforms of measures as some functionals of Fourier transforms, \emph{Probab. Math. Stat.} vol. \textbf{26} , Fasc. 1, pp. 187-200.

\medskip
Z. J.  Jurek (2007), Random integral representations for free-infinitely divisible and tempered
stable distributions, \emph{Stat. $\&$ Probab. Letters,} \textbf{77},
no 4, pp. 417-425.

\medskip
 Z. J.   Jurek (2011), The random integral representation conjecture: a quarter of a century later,
\emph{Lithuanian Math. Journal}, \textbf{51}, no 3, pp. 362-369,

\medskip
Z. J. \ Jurek (2012), Calculus on random integral mappings $I^{h,r}_{(a,b]}$ and their domains,
arXiv:0817121 [math PR]

\medskip
Z. J. \ Jurek and W. Vervaat (1983), An integral representation for selfdecomposable Banach space valued random variables, \emph{Z. Wahrsch. verw. Gebiete} \textbf{62}, pp. 246-262.

\medskip
M. M. Meerscheart and H.P. Scheffler (2001),\emph{Limit Distributions for Sums of Independent Random Vectors},
Wiley Series in Probability and Statistics, John Wiley $\&$ Sons, Inc. New York.

\medskip
A. Nica and R. Speicher (2006), \emph{Lectures on the combinatorics of free probability},
London Mathematical Society, Lecture Note series 335, Cambridge University Press.

\medskip
\noindent Zbigniew J.\ Jurek \hfil\break
Institute of Mathematics \hfil\break
University of Wroc\l aw \hfil\break
Pl.\ Grunwaldzki 2/4 \hfil\break
50-386 Wroc\l aw \hfil\break
Poland \hfil\break
{\it E-mail: \/} zjjurek@math.uni.wroc.pl ,  \ \  www.math.uni.wroc.pl/$\sim$zjjurek   \hfil\break

\end{document}